\documentclass{article}
\begin{document}
\centerline{\large\bf On the Number of Representations of Integers by} 
\centerline{\large\bf various Quadratic and Higher Forms}

\[
\]

\begin{quote}

\centerline{N.D. Bagis and M.L. Glasser}

\[
\]

\centerline{\bf Abstract}
We give formulas for the number of representations of non negative integers by various quadratic forms. We also give evaluations in the case of sum of two cubes (cubic case) and the quintic case, as well. We introduce a class of generalized triangular numbers and give several evaluations. Finally, we present a mean value asymptotic formula for the number of representations of an integer as sum of two squares known as the Gauss circle problem.        
\[
\]
\textbf{keywords}: \textrm{Quadratic Forms; Diophantine Equations; Sums of Squares; Asymptotics; Cubic Form; Quintic Form; Special Functions}

\end{quote}

\section{Introduction.}
  
The study of quadratic forms has been built up by many great mathematicians such as Euler, Gauss, Dirichlet, Liouville, Eisenstein, Glaisher, Ramanujan among others. This theory has applications to a wide number of areas in modern mathematics including Gauss' circle problem in higher dimensions, class number theory, algebraic geometry, elliptic and theta functions, the Fermat-Wiles theorem, Eisenstein series and many other (see [2-10]).\\ In this article using simple arguments we try to address the problem.\\  
\\  
We start with $K(x)$, the complete elliptic integral of the first kind, given by
\begin{equation}
K(x)=\int^{\pi/2}_{0}\frac{d\theta}{\sqrt{1-x^2\sin^2(\theta)}}=\frac{\pi}{2}{}_2F_1\left(\frac{1}{2},\frac{1}{2};1;x^2\right),
\end{equation}
where ${}_2F_1$ is Gauss hypergeometric function.\\
In terms of Weber's $\lambda(\tau)$-modular function (see [3],[4]) 
\begin{equation}
\lambda(\tau)=16q\prod^{\infty}_{n=1}\left(\frac{1+q^{2n}}{1+q^{2n-1}}\right)^8,
\end{equation}
where $q=e^{i\pi \tau}$, $Im(\tau)>0$,
 $\tau=\sqrt{-r}$, the singular modulus  $k=k_r$, $r>0$ is  
\begin{equation}
k^2_{r}=\lambda(\tau)=\left(\frac{\theta_2(q)}{\theta_3(q)}\right)^4\textrm{, }
\end{equation}
with 
\begin{equation}
\theta_2(q)=\sum^{\infty}_{n=-\infty}q^{(n+1/2)^2}\textrm{ and }\theta_3(q)=\sum^{\infty}_{n=-\infty}q^{n^2}\textrm{, }|q|<1
\end{equation} 
Also $k=k_r$, $0<k<1$ is the solution of the equation 
\begin{equation}
\frac{K(\sqrt{1-k_r^2})}{K(k_r)}=\sqrt{r}
\end{equation}    
As usual we set $K=K(k_r)$ (the complete elliptic integral at singular values) and $K'=K(k'_r)$, where $k'_r=\sqrt{1-k_r^2}$ is the complementary singular modulus. The  Fourier expansion for the Jacobi elliptic function $\textrm{dn}$ (see [3] p.51-53) is\\
\begin{equation}
\textrm{dn}(q,u)=\frac{\pi}{2K}+\frac{2\pi}{K}\sum^{\infty}_{n=1}\frac{q^n}{1+q^{2n}}\cos(2nz)
\end{equation}  
where $z=\left(\frac{\pi}{2K}\right)u$ lies in the strip $|Im(z)|<\frac{\pi}{2}Im(\tau)$, $\tau=i\frac{K'}{K}$.\\

A very interesting connection between number theory and the theory of elliptic functions stems from  Jacobi's famous  theorem\\
\\
\textbf{Theorem 1.} (Jacobi [3])\\
If $q=e^{-\pi\sqrt{r}}$, $r>0$, then
\begin{equation}
\theta_3(q)=\sum^{\infty}_{n=-\infty}q^{n^2}=\sqrt{\frac{2K}{\pi}}
\end{equation}
\\
This theorem plays a key role in the theory of elliptic functions  and we shall use it here in our investigation of quadratic forms of general type.
It is very easy to see  by setting $u=0$ in (6), using $\textrm{dn}(q,0)=1$ and then multiplying both sides of (6) by $2K/\pi$, that  
$$
\frac{2K}{\pi}=1+4\sum^{\infty}_{m=1}\frac{q^m}{1+q^{2m}}=1+4\sum^{\infty}_{m=1}q^m\sum^{\infty}_{l=0}(-1)^lq^{2ml}=
$$
$$
=1+4\sum^{\infty}_{m=1}\sum^{\infty}_{l=0}(-1)^lq^{(2l+1)m}.
$$
Writing $n=(2l+1)m$, $d=2l+1$, so $l=(d-1)/2$, if $d$ runs through the odd divisors of $n$ we have 
\begin{equation}
\frac{2K}{\pi}=1+4\sum^{\infty}_{n=1}\left[\sum_{d-odd, d|n}(-1)^{\frac{d-1}{2}}\right]q^n
\end{equation}
We define $\delta_0(n)=1$, if $n=0$ and $\delta_0(n)=4\sum_{d-odd, d|n}(-1)^{\frac{d-1}{2}}$, if $n\geq 1$.
If $r(n)$ denotes the number of representations of $n$ by the form  
$$
n=x^2+y^2\textrm{, }(x,y\in\bf Z\rm),
$$
then, if we consider the fact that
$$
\theta_3(q)^2=\sum^{\infty}_{n=-\infty}q^{n^2}\sum^{\infty}_{m=-\infty}q^{m^2}=\sum^{\infty}_{n,m=-\infty}q^{n^2+m^2}=\sum^{\infty}_{n=0}r(n)q^n
$$
and apply Jacobi's Theorem 1, we get \\
\\
\textbf{Theorem 2.} (Jacobi [8])\\
For $n=1,2,\ldots$ we have 
\begin{equation}
r(n)=4\sum_{d-odd,\textrm{ }d|n}(-1)^{\frac{d-1}{2}}
\end{equation}
and $r(0)=1$.\\ 

\section{Generalizations of Jacobi's two-square theorem}

Suppose we have two positive integers $A,B$, with $\gcd(A,B)=1$, and let $r_{A,B}(n)$ denote the number of representations of $n$ by the quadratic form
\begin{equation}
n=Ax^2+By^2
\end{equation}
Then 
$$
\theta_3\left(q^A\right)^2\theta_3\left(q^B\right)^2=\left(\sum^{\infty}_{n,m=-\infty}q^{An^2+Bm^2}\right)^2=\left(\sum^{\infty}_{n=0}r_{A,B}(n)q^n\right)^2
$$
But, also
$$
\theta_3\left(q^A\right)^2\theta_3\left(q^B\right)^2=\left(\sum^{\infty}_{n=0}r(n)q^{nA}\right)\left(\sum^{\infty}_{m=0}r(m)q^{mB}\right)=
$$
$$
=\sum^{\infty}_{n=0}\left(\sum_{kA+lB=n}r(k)r(l)\right)q^n
$$
The linear Diophantine equation $kA+lB=n$ has solutions for all $n$ since $\gcd(A,B)=1|n$.\\
We now introduce the transformation $T$, which assigns the Taylor coefficient $f_n$ of a function $f(q)$ to the Taylor coefficient $\left(\sqrt{f}\right)_n$   of its square root $\sqrt{f(q)}$, i.e. 
$$
\left(\sqrt{f}\right)_n=T(f_n).
$$

The  transform  $T$ can  be evaluated  by using Faa Di Bruno's Formula (see [1] p.823), which in this case is
\begin{equation}
T(f_n)=\sum_{m=0}^{n} h_m(f_0)\sum'\prod^{n}_{j=1}\frac{f_j^{a_j}}{a_j!}
\end{equation}
where the prime on the sum means that we sum over all non-negative integers $a_j$ such that $a_1+2a_2+3a_3+\ldots+na_n=n$ and $a_1+a_2+a_3+\ldots+a_n=m$. The function 
$h_m(x)=(-1)^{m}x^{1/2-m}\left(\frac{-1}{2}\right)_m$, $(a)_m=\frac{\Gamma\left(a+m\right)}{\Gamma(a)}$, $m=1,2,\ldots$.\\
\\
With the above notation we can proceed to\\ 
\\
\textbf{Proposition 1.}\\
Given two positive integers $A$,$B$ with $\textrm{gcd}(A,B)=1$ the number of the representations of $n=1,2,\ldots$ by the form $Ax^2+By^2$ is exactly 
\begin{equation}
r_{A,B}(n)=T\left(\sum_{kA+lB=n}r(k)r(l)\right)
\end{equation}
\\
Note that $r_{A,B}(0)$ is obviously 1. 
\\
\textbf{Proposition 2.}\\
Given two positive integers $A$,$B$ with $\textrm{gcd}(A,B)=1$, the number of  representations of $n=1,2,\ldots$ by the form $Ax^2+By^2$ is exactly 
\begin{equation}
r_{A,B}(n)=\left[\frac{1}{n!}\frac{d^{n}}{dq^n}\sqrt{\sum^{n}_{t=0}\left(\sum_{kA+lB=t}r(k)r(l)\right)q^t}\right]_{q=0}
\end{equation}
\\
In the same way as above we can prove\\
\\
\textbf{Theorem 3.}\\
If $A_1,A_2,\ldots,A_N$ are positive integers such that $\textrm{gcd}(A_1,A_2,\ldots,A_N)=1$, the number of the representations of $n=1,2,\ldots$ by the form $\sum^{N}_{k=1}A_kx_k^2$ is exactly 
\begin{equation}
r_{2}(N,n)=T\left(\sum_{k_1A_1+k_2A_2+\ldots+k_NA_N=n}r(k_1)r(k_2)\ldots r(k_N)\right)
\end{equation}
and
$$
r_2(N,n)=
$$
\begin{equation}
=\left[\frac{1}{n!}\frac{d^n}{dq^n}\sqrt{\sum^{n}_{t=0}\left(\sum_{k_1A_1+k_2A_2+\ldots+k_NA_N=t}r(k_1)r(k_2)\ldots r(k_N)\right)q^t}\right]_{q=0}
\end{equation}
\\
\textbf{Proposition 3.}\\
Consider the non-homogeneous quadratic form  
\begin{equation}
Ax^2+By^2+Cx+Dy+E
\end{equation}
with $A,B$ positive integers, $C,D,E$ general integers,   $gcd(A,B)=1$, and\\ $C\equiv0\textrm{ }\textrm{mod}(2A)$, $D\equiv0\textrm{ }\textrm{mod}(2B)$. Then $n$ has exactly 
$$
r_{A,B}\left(n+\frac{C^2}{4A}+\frac{D^2}{4B}-E\right)
$$
representations by (16).\\
\\
\textbf{Proof.}\\
Write $C=-2L_1A$ and $D=-2L_2B$. Then $n=Ax^2+By^2+Cx+Dy+E$ is equivalent to $n=A(x-L_1)^2+B(y-L_2)^2-AL_1^2-BL_2^2+E$ and the number of representations of $n$ by (16) is equal to the number of representation of  $n+AL_1^2+BL_2^2-E=n+\frac{C^2}{4A}+\frac{D^2}{4B}-E$, by $Ax^2+By^2$. qed\\
\\
\textbf{Application 1.}\\
Let $A,B,C,D$ be as in Proposition 3, then
$$
\sum^{\infty}_{n=-\infty}q^{An^2+Cn}\cdot\sum^{\infty}_{n=-\infty}q^{Bn^2+Dn}
=2\pi^{-1}q^{-n_0}K(k_r)\sqrt{m_{A,r}m_{B,r}}
$$
where $n_0=\frac{C^2}{4A}+\frac{D^2}{4B}$ and $q=e^{-\pi\sqrt{r}}$. The function $m_{n,r}=\frac{K(k_{n^2r})}{K(k_r)}$ is called a multiplier (see [4] pg.136) and takes algebraic values when $n$ is a positive integer and $r$ is rational.\\
\\
\textbf{Proof.}\\
From Proposition 3 we have
$$
\sum^{\infty}_{n=-\infty}q^{An^2+Cn}\cdot\sum^{\infty}_{n=-\infty}q^{Bn^2+Dn}=\sum^{\infty}_{n,m=-\infty}q^{An^2+Bm^2+Cn+Dm}=
$$ 
$$
=\sum^{\infty}_{n=0}r_{A,B}(n)q^{n-n_0}
=q^{-n_0}\sum^{\infty}_{n=0}r_{A,B}(n)q^n
=q^{-n_0}\vartheta_3(q^A)\vartheta_3(q^B)=
$$
$$
=q^{-n_0}\sqrt{\frac{2K(k_{A^2r})2K(k_{B^2r})}{\pi^2}}
$$
$$
=q^{-n_0}\frac{2K}{\pi}\sqrt{m_{A,r}m_{B,r}}\textrm{. qed }
$$
\\
\textbf{Application 2.}\\
The equation
\begin{equation}
k(Ax^2+By^2+Cx+Dy+E)+l=n
\end{equation}
has $r=r_{A,B}\left(\frac{n-l}{k}+\frac{C^2}{4A}+\frac{D^2}{4B}-E\right)$ solutions.\\In general if $P_N(x)=\sum^{N}_{k=0}a_kx^k$ is a polynomial with integer coefficients
and there exists exactly one integer $n'$ such that $P_N(n')=n$ then
\begin{equation}
P_N\left(Ax^2+By^2+Cx+Dy+E\right)=n
\end{equation}
has 
\begin{equation} r_{A,B}\left(n'+\frac{C^2}{4A}+\frac{D^2}{4B}-E\right)
\end{equation}
integer solutions, (including 0).\\Furthermore, if the equation $P_N(n')=n$ has integer solutions $n'=n'_1,n'_2,\ldots,n'_s$  with $s\leq N$, then the number of representations of $n$ by (18) will be
\begin{equation}
r=\sum_{i=1}^{s}r_{A,B}\left(n'_i+\frac{C^2}{4A}+\frac{D^2}{4B}-E\right).
\end{equation}
Any non-integer solution $n'$ to $P_N(n')=n$ leads to no  representation (18) and hence makes no contribution to the sum (20).\\ 
\\

Consider now the function $\sum^{\infty}_{n=0}q^{n^{\nu}}$, $\nu\in\bf N\rm$ and $\nu>2$. Then  
\begin{equation}
\left(\sum^{\infty}_{n=0}q^{n^{\nu}}\right)^2=\sum^{\infty}_{t=0}\left(\sum_{a^{\nu}+b^{\nu}=t}1\right)q^t.
\end{equation}   
Set
\begin{equation}
\textbf{1}_{\nu}(t):=T\left(\sum_{a^{\nu}+b^{\nu}=t}1\right).
\end{equation}
Then $P_{\nu}(n):=\textbf{1}_{\nu}(n)=1$, when $n$ is if the form $m^{\nu}$, ($m$ positive integer) and 0 otherwise. This leads to\\ 
\\ 
\textbf{Theorem 4.}\\
The number of representations of $n$ by $x^{\nu}+y^{\nu}$, where $x,y$ are non-negative integers, is
\begin{equation}
r_{\nu}(n)=\sum^{n}_{k=0}P_{\nu}(k)P_{\nu}(n-k).
\end{equation}
\\ 
\textbf{Proof.}\\
From (22) we get
\begin{equation}
\sum_{a^{\nu}+b^{\nu}=n}1=T^{(-1)}\left(\textbf{1}_{\nu}(n)\right)
\end{equation}
where $T^{(-1)}(f_n)$ is the $n$-th Taylor coefficient of $f^2$. Hence from  the Leibniz formula 
\begin{equation}
T^{(-1)}(f_n)=\sum_{k+l=n}f_kf_l
\end{equation}
Also, in the case where $n$ is the $\nu$-th power of a positive integer,  formula (23) gives $r_{\nu}(n^{\nu})=0$, (Fermat-Wiles theorem). qed\\
\\   
In general if $A(n)$ is a polynomial with positive integer coefficients then  the  equation $A(a)+A(b)=n$, where $a,b,n$ are non negative integers, has 
\begin{equation}
\sum_{A(a)+A(b)=n}1=\sum^{n}_{k=0}G_{A}(k)G_{A}(n-k)
\end{equation}
solutions. The function $G_A(n)$ is such that $G_A(n)=1$ if there exists a positive integer $m$  such that $n=A(m)$ and 0 otherwise.

\section{Representations by some  cubic and quintic forms}

In this section we give two formulas similar to  Jacobi's (Theorem 2) for the representation of a positive integer by a cubic and by a quintic form. The results of Section 2 can be generalized to higher order  terms under certain conditions. Historically, there are some  results  known regarding the cubic case. For example it is known that the Diophantine equation
\begin{equation}
ax^3-by^3=n
\end{equation}      
for $a,b,n$ integers, has  a finite number of solutions (see [6]).\\
From  the  Fermat-Wiles theorem it is known that
\begin{equation}
x^3+y^3=z^3
\end{equation}      
has only trivial solutions i.e. $\{x,0,x\}$ and $\{0,x,x\}$.\\  
Also a result of Euler states that  the equation
\begin{equation}
x^3+y^3=z^2
\end{equation}
admits  a parametric solution in integers (see [5] p.578-579).\\
We proceed by stating and proving \\
\\  
\textbf{Theorem 5.}\\
The number representations of $n$ by the form $x^3+y^3$  ($x,y$ non negative integers) is 
\begin{equation}
r_3(n)=\sum_{\scriptsize
\begin{array}{cc} 
	d|n\\
	d^3-4n=0
\end{array}\normalsize}1+2\cdot\sum_{\scriptsize
\begin{array}{cc} 
	d|n\\
	d^3-4n\neq0
\end{array}\normalsize
}S\left(\frac{-d^2+4\frac{n}{d}}{3}\right)
\end{equation} 
where $S(n)=1$ if $n$ is perfect square and 0 otherwise.\\ 
\\ 
\textbf{Proof.}\\
One has $x^3+y^3=(x+y)(x^2-xy+y^2)$ so if we set $u=x+y$ and $v=x^2-xy+y^2$,  $x,y$ are given by 
$$
x=\frac{1}{6}\left(3u-\sqrt{-3u^2+12v}\right)\textrm{, }y=\frac{1}{6}\left(3u+\sqrt{-3u^2+12v}\right)
$$
Hence we get the necessary and sufficient conditions for $u,v$ to determine  an integer $n=uv$ that can be expressed as the sum of two cubes.\\ 
\\
The quintic case is the similar the cubic. We have\\
\\
\textbf{Theorem 6.}\\
The  number of representations of $n$ by the form $x^5+y^5$  ($x,y$ non-negative integers) is 
$r_5(0)=1$ and if $n$ positive integer
$$
r_5(n)=-\sum_{\scriptsize
\begin{array}{cc} 
	d|n\\
	d^5-16n=0
\end{array}\normalsize}1+
$$
$$
+2\cdot\sum_{\scriptsize
\begin{array}{cc} 
	d|n\\
	d^5-16n\neq0\\	
\end{array}\normalsize}\textbf{X}_{N}\left
(\frac{5d-\sqrt{-25d^2+10\sqrt{5d^4+20\frac{n}{d}}}}{10}\right)S\left(5d^4+20\frac{n}{d}\right)\times 
$$
\begin{equation}
\times S\left(-25d^2+10\sqrt{5d^4+20\frac{n}{d}}\right)
\end{equation}
where $\bf X_N\rm$ is the characteristic function on the positive integers.

\section{Generalized Triangular Numbers}

We call  
\begin{equation}
t_m(n)=\frac{n^2+mn}{2}\textrm{, }n=0,1,2,\ldots\textrm{, with }m=0,1,2,\ldots 
\end{equation}
the $m$-triangular numbers.
We interested in the number of representations of a certain non-negative integers $n$ as the sum of $N$ in $m$-triangular numbers
\begin{equation}
n=\sum^{N}_{k=1}t_m(x_k)=\sum^{N}_{k=1}\frac{x_k^2+mx_k}{2}\textrm{, where }x_k\in\textbf{Z}\textrm{ and }m=0,1,2,\ldots 
\end{equation}  
The case of $m=1$, $N=2,3,4,\ldots$ is the well known representation of $n$ into simple triangular numbers (1-triangular numbers) and has been treated by many mathematicians (see [11]). The case $m=0$ is Jacobi's $N$-square theorem. At this point we drop the notation $r_{A,B}(n)$ we used  above and denote the number of representations of $n$ in (33) by $r_{m,N}(n)$. Also we denote $r(n)$ of (8) as $r_2(n)$, but the symbol $r_N(n)$ is left as in previous sections i.e is the number of representations of $n$ by   the diagonal form 
$
\sum^{N}_{k=1}x_k^2
$. 
Also, we recall the definition of certain theta functions studied by Ramanujan (see [14] pg.36):\\
\\
\textbf{Definition 1.}\\
If $|q|<1$, then
\begin{equation}
\phi(q):=\sum^{\infty}_{n=-\infty}q^{n^2}=\frac{(-q;q^2)_{\infty}(q^2;q^2)_{\infty}}{(q;q^2)_{\infty}(-q^2;q^2)_{\infty}}
\end{equation}
\begin{equation}
\psi(q):=\sum^{\infty}_{n=0}q^{n(n+1)/2}=\frac{(q^2;q^2)_{\infty}}{(q;q^2)_{\infty}}
\end{equation}
\begin{equation}
f(-q):=\sum^{\infty}_{n=-\infty}(-1)^nq^{n(3n-1)/2}=(q;q)_{\infty}
\end{equation}
where
\begin{equation}
(a;q)_{\infty}:=\prod^{\infty}_{n=0}(1-aq^n).
\end{equation}

Consider now the Jacobi triple product formula (see [15] pg.169-172):
\begin{equation}
\sum^{\infty}_{n=-\infty}q^{n^2+zn}=\prod^{\infty}_{n=0}(1-q^{2n+2})(1+q^{2n+1-z})(1+q^{2n+1+z})
\end{equation}
where $|q|<1$.\\
In case  $z=2p+1$, with $p$ a non-negative integer we get  
$$
\sum^{\infty}_{n=-\infty}q^{n^2+(2p+1)n}=f(-q^2)\prod^{\infty}_{n=0}(1+q^{2(n-p)})(1+q^{2(n+p)+2})=
$$
$$
=f(-q^2)\prod^{\infty}_{n=0}(1+q^{2(n-p)})(1+q^{2(n+p)+2})=
$$
$$
=f(-q^2)\prod^{\infty}_{n=0}(1+q^{2n})\prod^{p-1}_{n=0}(1+q^{2(n-p)})\frac{\prod^{\infty}_{n=0}(1+q^{2n+2})}{\prod^{p-1}_{n=0}(1+q^{2n+2})}=
$$
$$
=2q^{-p(p+1)}f(-q^2)(-q^2;q^2)^2_{\infty}.
$$
Since
$$
\prod^{p-1}_{n=0}\frac{1+q^{2(n-p)}}{1+q^{2n+2}}=q^{-p(p+1)}
$$
and
$$
(-q^2;q^2)_{\infty}=\frac{(q^4;q^4)_{\infty}}{(q^2;q^2)_{\infty}}
$$
we have\\
\\
\textbf{Proposition 4.}\\
If $|q|<1$ and $p=0,1,2,\ldots$, then 
\begin{equation}
\sum^{\infty}_{n=-\infty}q^{t_{2p+1}(n)}=2q^{-p(p+1)/2}\frac {f(-q^2)^2}{f(-q)}=2q^{-p(p+1)/2}\psi(q)
\end{equation}
\\
\textbf{Proof.}\\
The first equality follows from above discussion. For the second equality we have
$$
\frac{f(-q^2)^2}{f(-q)}=\frac{(q^2;q^2)^2_{\infty}}{(q;q)}=\frac{(q^2;q^2)}{(q;q^2)}=\psi(q),
$$
since 
\begin{equation}
(q;q^2)_{\infty}\cdot(-q;q)_{\infty}=1
\end{equation}
\\
\textbf{Proposition 5.}\\
If $|q|<1$ and $p\in\textbf{Z}$, then
\begin{equation}
\sum^{\infty}_{n=-\infty}q^{t_{2p}(n)}=q^{-p^2/2}\phi(q^{1/2})
\end{equation}
\\
\textbf{Proof.}\\
The proof is elementary since 
$$
\sum^{\infty}_{n=-\infty}q^{(n+p)^2}=\phi(q)
$$
when $|q|<1$ and $p\in\textbf{Z}$.\\
\\

Using the above Propositions we can generalize all the results in [11]. Sttarting from the $2p+1$-triangular numbers and Proposition 4, we immediately have\\ 
\\
\textbf{Theorem 7.}\\
If  $N\geq 2$ is an integer and $\delta_N(n)$ denotes the the number of representations of $n$ as $N$ 1-triangular numbers  
\begin{equation}
n=\sum^{N}_{k=1}\frac{x_k^2+x_k}{2}
\end{equation}
and if $r_{2p+1,N}(n)$ denotes the number of representations of the positive integer $n$ as $N$ $2p+1$-triangular numbers
\begin{equation}
n=\sum^{N}_{k=1}\frac{x_k^2+(2p+1)x_k}{2},
\end{equation}
then
\begin{equation}
r_{2p+1,N}(n)=\delta_N\left(n+\frac{N p(p+1)}{2}\right).
\end{equation}
\\

For example we have\\
\\
\textbf{Example 1.}\\
\textbf{i)} The number of  ways of representing  $n$ in the form
\begin{equation}
n=\frac{x^2+(2p+1)x}{2}+\frac{y^2+(2p+1)y}{2}
\end{equation}
is
\begin{equation}
r_{2p+1,2}(n)=4d_1(8(n+p^2+p)+2)-4d_3(8(n+p^2+p)+2) 
\end{equation} 
where 
\begin{equation}
d_a(n)=\sum_{\scriptsize
\begin{array}{cc} 
	d|n\\
	d\equiv a(4)
\end{array}\normalsize}1\textrm{, with }a=1,3
\end{equation}
\\
\textbf{ii)}\\
\begin{equation}
r_{2p+1,4}(n)=\sigma_1\left(2n+4p(p+1)+1\right)
\end{equation} 
where $\sigma_{\nu}(n)=\sum_{d|n}d^{\nu}$ is the divisor function.\\
\\

Continuing from Section 2 we get expressions for $r_{2p,N}(n)$.\\Since we know
\begin{equation}
r_N(n)=T\left(\sum_{k_1+k_2+\ldots+k_N=n}r(k_1)r(k_2)\ldots r(k_N)\right),
\end{equation}
we obtain  from Proposition 5 \\ 
\\
\textbf{Theorem 8.}\\
It is known that $r_N(2n)$ is the number of ways to represent the positive integer $n$ as the sum of  $N$ 0-triangular numbers  
\begin{equation}
n=\sum^{N}_{k=1}\frac{x_k^2}{2}.
\end{equation} 
The number of  representations of $n$ by
\begin{equation}
n=\sum^{N}_{k=1}\frac{x_k^2+2px_k}{2}
\end{equation}
is
\begin{equation}
r_{2p,N}(n)=r_N\left(2n+Np^2\right).
\end{equation}
\\

Continuing in this way, from Jacobi's two-square theorem we know that $$
r_{2}(n)=\sum_{d-odd,d|n}(-1)^{\frac{d-1}{2}},
$$ 
if $n=1,2,\ldots$ and $r_2(0)=1$. Combining the above results we obtain  next \\
\\
\textbf{Theorem 9.}\\
The number of representations $s_m(n)$ of $n$ in the form
\begin{equation}
\frac{x^2+mx}{2}+\frac{y^2+my}{2}
\end{equation}
is\\
\textbf{i)} If $m$ is even
\begin{equation}
s_m(n)=4\sum_{\scriptsize
\begin{array}{cc} 
	d|\left(n+\frac{m^2}{4}\right)\\
	d\equiv1(2)
\end{array}\normalsize}(-1)^{\frac{d-1}{2}}
\end{equation}
\textbf{ii)} If $m$ is odd
\begin{equation}
s_m(n)=4\sum_{\scriptsize
\begin{array}{cc} 
	d|(m^2+4n)\\
	d\equiv1(4)
\end{array}\normalsize}1-4\sum_{\scriptsize
\begin{array}{cc} 
	d|(m^2+4n)\\
	d\equiv3(4)
\end{array}\normalsize}1
\end{equation}
\\
\textbf{Theorem 10.}\\
The number of representations of $n$ as a sum of four $m$-triangular numbers 
\begin{equation}
\frac{x^2+mx}{2}+\frac{y^2+my}{2}+\frac{z^2+mz}{2}+\frac{w^2+mw}{2}
\end{equation}
is\\
\textbf{i)} If $m=2p$, $p=0,1,2,\ldots$,  
\begin{equation}
r_{2p,4}(n)=r_4(2n+4p^2),
\end{equation}
where (see [12]): 
\begin{equation}
r_4(n)=8\sum_{d|n}d\textrm{, if }n\textrm{ is odd} 
\end{equation}
and
\begin{equation}
r_4(n)=24\sum_{\scriptsize
\begin{array}{cc} 
	d|n\\
	d\equiv1(2)
\end{array}\normalsize}d\textrm{, if }n\textrm{ is even}  
\end{equation}
\textbf{ii)} If $m=2p+1$, $p=0,1,2,\ldots$,
\begin{equation} 
r_{2p+1,4}(n)=\sigma_1\left(2n+4p(p+1)+1\right)
\end{equation}
\\

From Theorem 10 we get the next\\
\\
\textbf{Theorem 11.}\\
For any given integer $m$, each non negative integer $n$ can represented as the sum of four $m$-triangular numbers.\\
\\
\textbf{Note.} Theorem 11 is a generalization of the   Lagrange's famous four square theorem.\\
\\ 
\textbf{Theorem 12.}\\
The number of representation of $n$ in the form
\begin{equation}
n=\frac{x^2+2px}{2}+\frac{y^2+2py}{2}+\frac{z^2+2pz}{2}
\end{equation}
is  
\begin{equation}
r_{2p,3}(n)=r_3(2n+3p^2)
\end{equation}
where
\begin{equation}
r_3(n)=\left\{\begin{array}{cc}
         24h(-n), \mbox{ } n\equiv3(8)\\
         12h(-4n), \mbox{ } n\equiv 1,2,5,6(8)\\
	       0, \mbox{  } n\equiv7(8)
\end{array}\right\}
\end{equation}
and where $h(n)$ is the class number of $n$.

\section{An Exponential Method}

\textbf{Proposition 6.}\\
In general, if $q=e^{-2x}$, $x>0$ and $X(n)$ is arithmetic function then
\begin{equation}
\sum^{\infty}_{n=1}X(n)\frac{n^2}{\sinh^2(nx)}=-\frac{d^2}{dx^2}\log\left(\prod^{\infty}_{n=1}\left(1-e^{-2nx}\right)^{X(n)}\right)
\end{equation}
\\
\textbf{Proof.}\\
See [18]. qed\\
\\

Let
\begin{equation} 
\chi_0(n)=\left\{
\begin{array}{cc}

	-2 \textrm{ if } n\equiv

	1(mod4)
	\\
	 
	 3 \textrm{ if } n\equiv 
	 
	 2(mod4)
\\	 
	 
	-2 \textrm{ if } n\equiv 
	
	3(mod4)
\\	
	
	 1 \textrm{ if } n\equiv 
	 
	 0(mod4).	 
\end{array}
\right\} 
\end{equation}
Then if $q=e^{-2x}$ we get (see [17],[18])
\begin{equation}
\sum^{\infty}_{n=1}\chi_0(n)\frac{n^2}{\sinh^2(n x)}=-\frac{d^2}{dx^2}\log\left(\sum^{\infty}_{n=-\infty}q^{n^2}\right)=-\frac{d^2}{dx^2}\log\left(\theta_3\left(e^{-2x}\right)\right)
\end{equation} 
and
\begin{equation}
\sum^{\infty}_{n=1}\frac{(-1)^nn^2}{\sinh^2(n x)}=-\frac{d^2}{dx^2}\log\left(\sum^{\infty}_{n=-\infty}q^{n(n+1)/2}\right)
\end{equation}  
A relationship between the hyperbolic sine function series and theta functions is:
\begin{equation} 
\makebox{If}\, \chi_{k,h}(n):=\left\{
\begin{array}{cc}
  
	1 \textrm{, if } n\equiv 
	 
	0,k+h,k-h(mod2k)
\\	
	
	0 \textrm{, otherwise },	 

\end{array}
\right\} 
\end{equation}
then
\begin{equation}
\sum^{\infty}_{n=1}\frac{\chi_{k,h}(n)n^2}{\sinh^2(n x)}=-\frac{d^2}{dx^2}\log\left(\sum^{\infty}_{n=-\infty}(-1)^nq^{kn^2+hn}\right)
\end{equation} 
when $k>h$, $k\in\bf N\rm$, $h\in\bf Z\rm$.\\
Assume that not both $k,h$ are even or odd, then from [18] 
\begin{equation}
\log\left(\sum^{\infty}_{n=-\infty}(-1)^nq^{kn^2+hn}\right)=-\sum^{\infty}_{n=1}f_{k,h}(n)q^n
\end{equation}
where 
\begin{equation}
f_{k,h}(n):=\frac{1}{n}\sum_{d|n}\chi_{k,h}(d)d.
\end{equation}
Hence 
\begin{equation}
\sum^{\infty}_{n=-\infty}q^{kn^2+hn}=\exp\left(- \sum^{\infty}_{n=1}(-1)^nf_{k,h}(n)q^n\right).
\end{equation}
If we define $T_0(a_n)$ to be such that
\begin{equation}
\exp\left(-\sum^{\infty}_{n=1}a_nx^n\right)=\sum^{\infty}_{n=0}T_0\left(a_n\right)x^n,
\end{equation}
then\\ 
\\
\textbf{Theorem 13.}\\
The number of representations of $n\in\bf N\rm$ in the form  
\begin{equation}
n=\sum^{N}_{l=1}k_lx_l^2+h_lx_l
\end{equation}
with $k_l>|h_l|>0$, $k_l,h_l$ not both even or both odd, $\forall l=1,2,\ldots,N$ is
\begin{equation}
r(n)=T_0\left(\sum^{N}_{l=1}(-1)^lf_{k_l,h_l}(n)\right)
\end{equation} 
where $f_{k,h}$ is that of (71) and $\chi_{k,h}$ is that of (68).\\ 
\\
\textbf{Examples.}\\
\textbf{i)} For example, the number of representations of $n$ in the form  
\begin{equation}
10x^2+11y^2+x+4y
\end{equation} 
is 
\begin{equation}
r(n)=T_0\left(\frac{(-1)^n}{n}\sum_{d|n}\chi_{10,1}(d)d+\frac{(-1)^n}{n}\sum_{d|n}\chi_{11,4}(d)d\right)
\end{equation}
\\
\textbf{ii)} Another example is the number of representations of $n$ by the form  
\begin{equation}
3x^2-2x+3y^2-2y
\end{equation}   
which is
\begin{equation}
r(n)=T_0\left(2(-1)^n\sigma^{*}_{6}(n)+2(-1)^n I_{6}(n)\right),
\end{equation}
with
\begin{equation} \sigma^{*}_{a}(n)=\frac{1}{n}\sum_{d|n}\left(\frac{d}{a^2}\right)d
\end{equation}
where $I_{a}(ka)=1$ for $k=1,2,\ldots$ and  is otherwise  0.\\
\\
\textbf{iii)} For
\begin{equation}
n=2x^2-x+2y^2-y
\end{equation}
we have
\begin{equation}
r(n)=T_0\left(2(-1)^n\sigma^{*}_{2}(n)+2(-1)^nI_{4}(n)\right)
\end{equation}
\\
\textbf{iv)} For
\begin{equation}
n=4x^2+3y^2+3x+2y
\end{equation}
we have
\begin{equation}
r(n)=T_0\left(\frac{(-1)^n}{n}\sum_{d|n}\chi_{4,3}(d)d+\frac{(-1)^n}{n}\sum_{d|n}\chi_{3,2}(d)d\right).
\end{equation}

\section{Asymptotic Expansion of $\sum_{n\leq x}r_2(n)$}

In this section we provide asymptotic formulas relating to the mean value of $r_2(n)$, using a formula of Hardy (see [16]). 
\begin{equation}
\sum_{n\leq x}r_2(n)=\pi x+x^{1/2}\sum^{\infty}_{n=1}\frac{r_2(n)}{\sqrt{n}}J_1(2\pi \sqrt{nx})
\end{equation}

If $a,b\in\textbf{R}$, then we define
\begin{equation}
M_s(a,b)=\sum^{\infty}_{k-odd\textrm{, }k=1}(-1)^{\frac{k+1}{2}}\frac{\cos(a+b\sqrt{k})}{k^s}
\end{equation}
\begin{equation}
N_s(a,b)=\sum^{\infty}_{k-odd\textrm{, }k=1}(-1)^{\frac{k+1}{2}}\frac{\sin(a+b\sqrt{k})}{k^s}
\end{equation}
and
\begin{equation}
P_{s}(a,b)=\sum^{\infty}_{n=1}\frac{M_{s}(a,b\sqrt{n})}{n^s}\textrm{, }Q_{s}(a,b)=\sum^{\infty}_{n=1}\frac{N_{s}(a,b\sqrt{n})}{n^s}.
\end{equation}
Next we prove \\ 
\\
\textbf{Theorem 14.}
$$
R(x)=\sum_{n\leq x}r_2(n)-x\pi=\frac{x^{1/4}}{\pi}P_{3/4}\left(\frac{\pi}{4},2\pi\sqrt{x}\right)+\sum^{N}_{s=1}\frac{(-1)^sc_1(2s)P_{s+3/4}\left(\frac{\pi}{4},2\pi\sqrt{x}\right)}{2^{4s}\pi^{2s+1}x^{s-1/4}}-
$$
\begin{equation}
-\sum^{N}_{s=0}\frac{(-1)^sc_1(2s+1)Q_{s+5/4}\left(\frac{\pi}{4},2\pi\sqrt{x}\right)}{2^{4s+2}\pi^{2s+2}x^{s+1/4}}+O\left(c_1(2N)4^{-N}x^{-N-1/2}\right)
\end{equation}
where $c_1(m)=(-1)^m\frac{\left(-\frac{1}{2}\right)_m\left(\frac{3}{2}\right)_m}{m!}$.\\
\\
\textbf{Proof.}\\
From (85) and (9) we have
$$
\sqrt{x}\sum_{n=1}^{\infty}\frac{r_2(n)}{\sqrt{n}}J_1(2\pi\sqrt{nx})=
\sqrt{x}\sum^{\infty}_{n=1}\left(\sum_{d-odd, d|n}(-1)^{\frac{d-1}{2}}\right)\frac{1}{\sqrt{n}}J_1(2\pi\sqrt{nx})=
$$
$$
=\sqrt{x}\sum^{\infty}_{n=1}\sum^{\infty}_{m=1}\frac{(-1)^m}{\sqrt{n(2m-1)}}J_1\left(2\pi\sqrt{n(2m-1)x}\right)=
$$
\begin{equation}
=\sqrt{x}\sum^{\infty}_{p-odd, n, p=1}\frac{(-1)^{\frac{p-1}{2}}}{\sqrt{np}}J_1\left(2\pi\sqrt{npx}\right)
\end{equation}
The function $J_1(x)$ has the following asymptotic expansion as $x\rightarrow\infty$
$$
J_1(x)=\sqrt{\frac{2}{\pi x}}[\cos\left(x-\frac{3\pi}{4}\right)\sum^{\infty}_{n=0}\frac{(-1)^nc_1(2n)}{(2x)^{2n}}-
$$
\begin{equation}
-\sin\left(x-\frac{3\pi}{4}\right)\sum_{n=0}^{\infty}\frac{(-1)^nc_1(2n+1)}{(2x)^{2n+1}}]
\end{equation}
The error due to stopping the summation at any term is the order of  magnitude of that term multiplied by $1/x$. Hence, using (91) in (90) we get (89).\\
\\

Setting $N=1$ in (89), we get
$$
R(x)=-\sum^{\infty}_{\scriptsize
\begin{array}{cc} 
	n,p=1\\
	p-odd
\end{array}\normalsize}[\frac{105 (-1)^{\frac{p+1}{2}}\sin\left(2\pi  \sqrt{npx}+\frac{\pi}{4}\right)}{4096\pi^3 (np)^{9/4} x^{5/4}}-\frac{15(-1)^{\frac{p+1}{2}} \cos \left(2\pi  \sqrt{npx}+\frac{\pi}{4}\right)}{256\pi^2 (np)^{7/4} x^{3/4}}+
$$
$$
+\frac{3(-1)^{\frac{p+1}{2}} \sin\left(2\pi  \sqrt{npx}+\frac{\pi}{4}\right)}{8\pi (np)^{5/4} \sqrt[4]{x}}-\frac{2 (-1)^{\frac{p+1}{2}} \sqrt[4]{x} \cos \left(2\pi\sqrt{np x}+\frac{\pi }{4}\right)}{(np)^{3/4}}]+O\left(x^{-3/4}\right)
$$
where $p=2l+1$. Therefore, we have\\
\\
\textbf{Proposition 7.}\\
The Gauss circle problem reduces to finding the rate of convergence of $R(x)=\frac{1}{x^{1/4}}\left(\sum_{n\leq x}r_{2}(n)-\pi x\right)$, which is equivalent to that of
\begin{equation}
S(x)=\sum^{\infty}_{n=1}\frac{r_2(n)\cos\left(2\pi\sqrt{nx}+\frac{\pi}{4}\right)}{n^{3/4}}=\sum^{\infty}_{n,l=1}\frac{(-1)^{l-1}\cos \left(2\pi\sqrt{n(2l-1)x}+\frac{\pi}{4}\right)}{(n(2l-1))^{3/4}}.
\end{equation} 
\\

We can write
$$
S(x)=\sum^{\infty}_{n,l=1}
\frac{1}{n^{3/4}}\left(\frac{\cos\left(2\pi\sqrt{n(4l+1)x}+\frac{\pi}{4}\right)}{(4l+1)^{3/4}}-\frac{\cos\left(2\pi\sqrt{n(4l-1)x}+\frac{\pi}{4}\right)}{(4l-1)^{3/4}}\right)
$$
Also if we set
\begin{equation}
\theta(n,l,x)=\cos\left(2\pi\sqrt{n(4l-1)x}+\frac{\pi}{4}\right)-\cos\left(2\pi\sqrt{n(4l+1)x}+\frac{\pi}{4}\right)
\end{equation}
then,  because
$$
\lim_{l\rightarrow\infty}(4l+1)^{11/4+1}\left(\frac{1}{(4l-1)^{3/4}}-\frac{1}{(4l+1)^{3/4}}-\frac{3/2}{(4l+1)^{3/4+1}}\right)=\frac{21}{8},
$$
we can write
\begin{equation}
\frac{1}{(4l-1)^{3/4}}-\frac{1}{(4l+1)^{3/4}}=\frac{3/2}{(4l+1)^{3/4+1}}+O\left(\frac{1}{(4l+1)^{11/4}}\right).
\end{equation}
Hence we get 
$$
S(x)=\sum^{\infty}_{n=1}\sum^{\infty}_{l=1}\frac{\cos\left(2\pi\sqrt{n(4l+1)x}+\frac{\pi}{4}\right)}{n^{3/4}}\left[\frac{1}{(4l+1)^{3/4}}-\frac{1}{(4l-1)^{3/4}}\right]+
$$
$$
+\sum^{\infty}_{n,l=1}\frac{\theta(n,l,x)}{n^{3/4}(4l+1)^{3/4}}.
$$
But, if we set $\sigma'_{\nu}(n)=\sum_{d\equiv1(4),d|n}d^{\nu}$ then (since the sequences we use involving the big $O$-symbol are bounded above; also see  [9] pg.135-136):
$$ \sum^{\infty}_{n,l=1}\frac{\cos\left(2\pi\sqrt{n(4l+1)x}+\frac{\pi}{4}\right)}{n^{3/4}}\left[\frac{1}{(4l-1)^{3/4}}-\frac{1}{(4l+1)^{3/4}}\right]=
$$
$$
\frac{3}{2}\sum^{\infty}_{n=1}\left[\sum^{\infty}_{l=1}\frac{\cos\left(2\pi
\sqrt{n(4l+1)x}+\frac{\pi}{4}\right)}{n^{3/4}(4l+1)^{3/4+1}}+O\left(\sum^{\infty}_{l=1}\frac{\cos\left(2\pi\sqrt{n(4l+1)x}+\frac{\pi}{4}\right)}{n^{3/4}(4l+1)^{11/4}}\right)\right]=
$$
$$
\frac{3}{2}\sum^{\infty}_{n=1}\sum^{\infty}_{l=1}\frac{n\cos\left(2\pi
\sqrt{n(4l+1)x}+\frac{\pi}{4}\right)}{n^{3/4+1}(4l+1)^{3/4+1}}+\sum^{\infty}_{n=1}O\left(\sum^{\infty}_{l=1}\frac{\cos\left(2\pi\sqrt{n(4l+1)x}+\frac{\pi}{4}\right)}{n^{3/4}(4l+1)^{11/4}}\right)=
$$
$$
\frac{3}{2}\sum^{\infty}_{n=1}\frac{\sigma'_1(n)\cos\left(2\pi
\sqrt{nx}+\frac{\pi}{4}\right)}{n^{3/4+1}}+
O\left( \sum^{\infty}_{n=1}\frac{\sigma'_2(n)\cos\left(2\pi\sqrt{nx}+\frac{\pi}{4}\right)}{n^{11/4}}\right)=
$$
\begin{equation}
\frac{3}{2}\sum^{\infty}_{n=1}\frac{\sigma'_1(n)}{n^{1+\delta}}\frac{\cos\left(2\pi
\sqrt{nx}+\frac{\pi}{4}\right)}{n^{3/4-\delta}}+O\left( \sum^{\infty}_{n=1}\frac{\sigma'_2(n)}{n^{2+\delta}}\frac{\cos\left(2\pi\sqrt{nx}+\frac{\pi}{4}\right)}{n^{3/4-\delta}}\right)
\end{equation}
Now if we keep in mind the inequality $\sigma'_{a}(n)=o(n^{a+\delta})$ and assume that the sums 
\begin{equation}
D_M(x):=\sum^{M}_{n=1}\frac{\cos\left(2\pi\sqrt{nx}+\frac{\pi}{4}\right)}{n^{3/4-\delta}} 
\end{equation}
are uniformly bounded, when $M\rightarrow\infty$, for every $\delta>0$ sufficiently small, then using Abel's test (see [13] pg.346) the two series in (95) are uniformly convergent.\\ 
Also  
$$
\left|\sum^{\infty}_{n,l=1}\frac{\theta(n,l,x)}{n^{3/4}(4l+1)^{3/4}}\right|\leq
$$
$$ \left|\sum^{\infty}_{n,l=1}\frac{\cos\left(2\pi\sqrt{n(4l+1)x}+\frac{\pi}{4}\right)}{(n(4l+1))^{3/4}}\right|+C\left|\sum^{\infty}_{n,l=1}\frac{\cos\left(2\pi\sqrt{n(4l-1)x}+\frac{\pi}{4}\right)}{(n(4l-1))^{3/4}}\right|=
$$
$$
=O\left(\sum^{\infty}_{n,l=1}\frac{\cos\left(2\pi\sqrt{n(4l+1)x}+\frac{\pi}{4}\right)}{(n(4l+1))^{3/4}}\right)
$$
\begin{equation}
=O\left(\sum^{\infty}_{n=1}\frac{\sigma'_{0}(n)}{n^{\delta}}\frac{\cos\left(2\pi\sqrt{nx}+\frac{\pi}{4}\right)}{n^{3/4-\delta}}\right)= O\left(\sum^{\infty}_{n=1}\frac{1}{n^{\delta-\epsilon}}\frac{\cos\left(2\pi\sqrt{nx}+\frac{\pi}{4}\right)}{n^{3/4-\delta}}\right)
\end{equation}
Since it is known that exists $C_0$ such that $\sigma'_{0}(n)\leq C_0n^{\epsilon}$, for $0<\epsilon<\delta$, if the sum (96) is  uniformly bounded,  again  using Abel's test we get the uniform convergence of (97).\\
\\ 
The Euler-Maclaurin formula for a function $F$ having 4 continuous derivatives in the interval $(a,m)$ states
$$
\sum^{M}_{k=a}F(k)=\int^{M+1}_{a}F(t)dt+\frac{1}{2}\left(F(M+a)+F(a)\right)+
$$
\begin{equation}
+\frac{1}{12}(F'(M+a)-F'(a))-\frac{1}{120}\sum^{M-1}_{k=0}F^{(4)}(a+k+\xi)
\end{equation}  
with $0<\xi<1$.\\
If we set 
$$
F(t)=\frac{\cos\left(2\pi\sqrt{tx}+\frac{\pi}{4}\right)}{\sqrt{t}}
$$ 
then
$$
\int^{M+1}_{a}\frac{\cos\left(2\pi\sqrt{(t-1)x}+\frac{\pi}{4}\right)}{\sqrt{t-1}}dt=\frac{\sin\left(2\pi\sqrt{xM}+\frac{\pi}{4}\right)}{\pi\sqrt{x}}-\frac{\sin\left(2\pi\sqrt{x}+\frac{\pi}{4}\right)}{\pi\sqrt{x}}.
$$
Also
$$
F'(t)=-\frac{\cos\left(2\pi\sqrt{x t}+\frac{\pi}{4}\right)}{2 t^{3/2}}-\sqrt{x}\pi\frac{\sin\left(2\pi\sqrt{xt}+\frac{\pi}{4}\right)}{t}
$$
and
$$
F^{(4)}(t)=\frac{\pi^4 x^2 \cos \left(\frac{1}{4} \pi  \left(8 \sqrt{t x}+1\right)\right)}{t^{5/2}}-\frac{5 \pi ^3 x \sqrt{t x} \sin \left(\frac{1}{4} \pi  \left(8 \sqrt{t x}+1\right)\right)}{t^{7/2}}+
$$
$$
+\frac{105 \pi  \sqrt{t x} \sin \left(\frac{1}{4} \pi  \left(8 \sqrt{t x}+1\right)\right)}{8 t^{9/2}}-\frac{45 \pi ^2 x \cos \left(\frac{1}{4} \pi  \left(8 \sqrt{t x}+1\right)\right)}{4 t^{7/2}}+
$$
$$
+\frac{105 \cos \left(\frac{1}{4} \pi  \left(8 \sqrt{t x}+1\right)\right)}{16 t^{9/2}}.
$$
Hence the Euler-Maclaurin summation formula assures us that the $\lim_{M\rightarrow\infty}D_M(x)$ exists i.e. the series $D(x)$ converges. Note that we don't use $n^{3/4}$ in the sum $D_M(x)$, but $\sqrt{n}$. One can see that this don't changes nothing (we proceed with $\sqrt{n}$ instead of $n^{3/4}$ for avoid showing large formulas).\\ 
\\
From all the  arguments in the present paragraph we are  able to prove \\
\\
\textbf{Proposition 8.}\\
For $\delta>0$ small enough the series
\begin{equation}
D(x)=\sum^{\infty}_{n=1}\frac{\cos\left(2\pi\sqrt{nx}+\frac{\pi}{4}\right)}{n^{3/4-\delta}} 
\end{equation} 
is convergent.\\ 
Further if exist always fixed $\delta_0>0$ small enough, such  for every $0<\delta<\delta_0$ the partial sums
\begin{equation}
D_M(x)=\sum^{M}_{n=1}\frac{\cos\left(2\pi\sqrt{nx}+\frac{\pi}{4}\right)}{n^{3/4-\delta}} 
\end{equation} 
are uniformly bounded, then
for every function $f(x)$ such that $\lim_{x\rightarrow \infty}f(x)=+\infty$ holds 
\begin{equation}
\sum_{n\leq x}r_2(n)=x\pi+O\left(x^{1/4}f(x)\right)\textrm{,  as }x\rightarrow\infty
\end{equation}
\\
\textbf{Proof.}\\
If $D_M(x)$ is uniformly bounded then it is also  uniformly convergent (because of parameter $\delta$) and we have
$$
\lim_{x\rightarrow\infty}\frac{D(x)}{f(x)}=\lim_{x\rightarrow\infty}\sum^{\infty}_{n=1}\frac{\cos\left(2\pi\sqrt{nx}+\frac{\pi}{4}\right)}{n^{3/4-\delta}f(x)}=\sum^{\infty}_{n=1}\lim_{x\rightarrow\infty}\frac{\cos\left(2\pi\sqrt{nx}+\frac{\pi}{4}\right)}{n^{3/4-\delta}f(x)}=0.
$$  
An example of such function $f(x)$ is $\log_n(x)=\underbrace{\log(\log(...(\log(}_{n-times}x))))$, for fixed large $n$.\\
\\
\textbf{Lemma 1.}(see [9] pg.145)\\
Suppose that $\lambda_1,\lambda_2,\ldots$ is a nondecreasing sequence of real numbers with limit infinity, that $c_1,c_2,\ldots$ is an arbitrary sequence of real or complex numbers, and that $f(x)$ has a continuous derivative for $x\geq \lambda_1$. Put
\begin{equation}
C(x)=\sum_{\lambda_n\leq x}c_n
\end{equation}   
where the summation is over all $n$ for which $\lambda_n\leq x$. Then for $x\geq \lambda_1$, 
\begin{equation}
\sum_{\lambda_n\leq x}c_nf(\lambda_n)=C(x)f(x)-\int^{x}_{\lambda_1}C(t) f'(t)dt
\end{equation}
\\
\textbf{Theorem 15.}\\
For every $\delta>0$ sufficiently small the sum (100) is always bounded (uniformly bounded) in $M$ and $x$. By this the asymptotic formula (101) is true and the Gauss circle problem is solved.\\
\\
\textbf{Proof.}\\
Let the function
\begin{equation}
f(t)=\frac{\cos\left(2\pi t\sqrt{a}+\frac{\pi}{4}\right)}{t^{3/2}}
\end{equation} 
It's first derivative is
\begin{equation}
f'(t)=-\frac{2 \pi  \sqrt{a} \sin \left(2 \pi  \sqrt{a} t+\frac{\pi }{4}\right)}{t^{3/2}}-\frac{3 \cos \left(2 \pi  \sqrt{a} t+\frac{\pi }{4}\right)}{2 t^{5/2}}
\end{equation}
From Lemma seting $y=\sqrt{M}$ we have
\begin{equation}
\sum_{\sqrt{n}\leq y}1=y^2=M
\end{equation}
Also with $y=\sqrt{M}$ we have
$$
\sum_{\sqrt{n}\leq y}f\left(\sqrt{n}\right)=\left(\sum_{\sqrt{n}\leq x}1\right)f(x)-\int^{x}_{1}t^2 f'(t)dt=
$$
$$
=\sum_{n\leq M}f\left(\sqrt{n}\right)=Mf(\sqrt{M})-\int^{\sqrt{M}}_{1}t^2 f'(t)dt
$$
Hence
$$
\sum_{n\leq M}f\left(\sqrt{n}\right)=\sum_{n\leq M}\frac{\cos\left(2\pi\sqrt{n a}+\frac{\pi}{4}\right)}{n^{3/4}}
=Mf\left(\sqrt{M}\right)-\int^{\sqrt{M}}_{1}t^2f'(t)dt=
$$
$$
=\frac{1}{\sqrt{2} \sqrt[4]{a}}\left(-2 F_C\left(2 \sqrt[4]{a}\right)+2 F_C\left(2 \sqrt[4]{a M} \right)+2 F_S\left(2 \sqrt[4]{a}\right)-2 F_S\left(2 \sqrt[4]{aM}\right)\right)+
$$
\begin{equation}
+\frac{1}{\sqrt{2}}\left(\cos \left(2 \pi  \sqrt{a}\right)-\sin \left(2 \pi  \sqrt{a}\right)\right)
\end{equation}
where $F_C(z)=\int^{z}_{0}\cos\left(\frac{\pi t^2}{2}\right)dt$ and $F_S(z)=\int^{z}_{0}\sin\left(\frac{\pi t^2}{2}\right)dt$ are the Fresnel$-C,S$ functions.\\
But function (107) is absolutely bounded when $M=1,2,\ldots$ and $a>0$ by some universal constant (we mean $\sqrt{2}$). Hence for the sum 
\begin{equation}
G(h,x,M):=\sum^{M}_{n=1}\frac{\cos\left(2\pi\sqrt{nx}+\frac{\pi}{4}\right)}{n^{3/4-h}}
\end{equation}
it holds that $G\left(0,x,M\right)$ is bounded in $M\in\textbf{N}$ and $x>0$. Using the mean-value theorem there exists $0<\xi<\delta$ such that
$$
\left|\frac{G(\delta,x,M)-G\left(0,x,M\right)}{\delta}\right|=\left|\partial_{h}G(\xi,x,M)\right|
$$
Hence
$$
|G(\delta,x,M)|\leq\delta\left|\partial_{h}G(\xi,x,M)\right|+\left|G\left(0,x,M\right)\right|=
$$
\begin{equation}
=\delta\left|\sum^{M}_{n=1}\frac{\cos\left(2\pi\sqrt{nx}+\frac{\pi}{4}\right)}{n^{1/2}}\frac{\log(n)}{n^{1/4-\xi}}\right|+\left|G\left(0,x,M\right)\right|
\end{equation}
If we consider next the function
\begin{equation}
f_1(t)=\frac{\cos\left(2\pi t\sqrt{a}+\frac{\pi}{4}\right)}{t}
\end{equation}
we can show as above that
\begin{equation}
\sum_{n\leq M}f_1(\sqrt{n})
\end{equation}   
is uniformly bounded and (109) becomes
$$
\left|G\left(\delta,x,M\right)\right|\leq \delta\left|\sum^{M}_{n=1}\frac{\cos\left(2\pi\sqrt{nx}+\frac{\pi}{4}\right)}{n^{1/2}}\frac{\log(n)}{n^{1/4-\xi}}\right|+\left|G\left(0,x,M\right)\right|=
$$
\begin{equation}
=O\left(\delta\left|\sum^{M}_{n= 1}\frac{\cos\left(2\pi\sqrt{nx}+\frac{\pi}{4}\right)}{n^{1/2}}\right|\right)+|G(0,x,M)|<\infty
\end{equation}
uniformly in $M$ and $x$, when $0<\delta<\frac{1}{4}$.\\
Hence $D_{M}(x)$ are bounded in $M$ and $x$ when $\delta$ is sufficiently small, and the theorem is proved.

\[
\]

\centerline{\bf References}\vskip .2in

\noindent

[1]: M. Abramowitz and I.A. Stegun. 'Handbook of Mathematical Functions'. Dover Publications, New York. (1972)

[2]: T. Apostol. 'Introduction to Analytic Number Theory'. Springer Verlag, New York, Berlin, Heidelberg, Tokyo, (1974)

[3]: J.V. Armitage W.F. Eberlein. 'Elliptic Functions'. Cambridge University Press. (2006)

[4]: J.M. Borwein and P.B. Borwein. 'Pi and the AGM'. John Wiley and Sons, Inc. New York, Chichester, Brisbane, Toronto, Singapore. (1987)

[5]: L.E. Dickson. 'History of the Theory of Numbers, Vol2: Diophantine Analysis'. Dover. New York, (2005)

[6]: G.H. Hardy. 'Ramanujan Twelve Lectures on Subjects Suggested by his Life and Work, 3rd ed.' Chelsea. New York, (1999)

[7]: M.D. Hirschhorn. 'Three classical results on representations of a number'. Seminaire Lotharingien de Combinatoire, (42) (1999)   

[8]: C.G.J. Jacobi. 'Fundamenta Nova Functionum Ellipticarum'. Werke I, 49-239. (1829) 

[9]: William J. LeVeque. 'Fundamentals of Number Theory'. Dover Publications. New York. (1996)

[10]: E.T. Whittaker and G.N. Watson. 'A course on Modern Analysis'. Cambridge U.P. (1927)

[11]: Ken Ono. 'Representations of Integers as Sums of Squares'.Journal of Number Theory. (95), 253-258. (2002)

[12]: Ila Varma. 'Sums of Squares, Modular Forms, and Hecke Characters'. Master thesis. Mathematisch Instituut, Universiteit Leiden. June 18 (2010). 

[13]: Konrad Knopp. 'Theory and Applications of Infinite Series'. Dover Publications, Inc. New York. (1990).

[14]: Bruce C. Berndt. 'Ramanujan`s Notebooks Part III'. Springer Verlag, New York (1991).

[15]: G.E. Andrews, Number Theory. Dover Publications, New York, (1994). 

[16]: G.H. Hardy. 'On the expression of a number as the sum of two squares'. Quart. J. Math. (Oxford) 46 (1915),
263–283.

[17]: N. Bagis. 'Some New Results on Sums of Primes'. Mathematical Notes, (2011), Vol. 90, No. 1, pp 10-19.

[18]: N. Bagis. 'Some Results on Infinite Series and Divisor Sums'.\\ arXiv:0912.48152v2 [math.GM] (2014).

\end{document}